\documentclass[12pt]{article}
\usepackage{amsfonts}
\usepackage{amssymb}
\usepackage{amsfonts,amssymb, amsmath, latexsym}
\usepackage{pst-all}
\usepackage{CJK}
\usepackage{mathrsfs}

\usepackage[left,pagewise,mathlines]{lineno}

\newtheorem{cor}{Corollary}[section]

\newtheorem{thm}{Theorem}[section]

\newtheorem{obs}[thm]{Observation}

\newenvironment{pf}[1][Proof]{\noindent\textbf{#1.} }{\hfill\rule{1mm}{2mm}}

\makeatletter \@addtoreset{equation}{section} \makeatother

\def\g{\gamma}

\def\br{\bar}

\begin{document}

\title{On the $p$-reinforcement and the complexity\thanks{The work was
supported by NNSF of China (No.10711233) and the Fundamental
Research Fund of NPU (No. JC201150)}}

\author{
{ You Lu$^a$ \quad Fu-Tao Hu$^b$} \quad Jun-Ming Xu$^b$\footnote{Corresponding author:
xujm@ustc.edu.cn (J.-M. Xu)} \\
{\small $^a$Department of Applied Mathematics,}\\
{\small             Northwestern Polytechnical University,}\\
{\small             Xi'an Shanxi 710072, P. R. China}\\
{\small             Email: luyou@nwpu.edu.cn}\\ \\
{\small $^b$Department of Mathematics,}\\
{\small             University of Science and Technology of China,}   \\
{\small             Wentsun Wu Key Laboratory of CAS,}\\
{\small             Hefei, Anhui, 230026, P. R. China}\\
{\small             Email: hufu@mail.ustc.edu.cn;\ xujm@ustc.edu.cn}  \\
}
\date{}
\maketitle

\begin{abstract}
Let $G=(V,E)$ be a graph and $p$ be a positive integer. A subset
$S\subseteq V$ is called a $p$-dominating set if each vertex not in
$S$ has at least $p$ neighbors in $S$. The $p$-domination number
$\g_p(G)$ is the size of a smallest $p$-dominating set of $G$. The
$p$-reinforcement number $r_p(G)$ is the smallest number of edges
whose addition to $G$ results in a graph $G'$ with $\g_p(G')<
\g_p(G)$. In this paper, we give an original study on the
$p$-reinforcement, determine $r_p(G)$ for some graphs such as paths,
cycles and complete $t$-partite graphs, and establish some upper
bounds of $r_p(G)$. In particular, we show that the decision problem
on $r_p(G)$ is NP-hard for a general graph $G$ and a fixed integer
$p\geq 2$.
\end{abstract}



\newpage

\section{Induction}

For notation and graph-theoretical terminology not defined here we
follow \cite{x03}. Specifically, let $G=(V,E)$ be an undirected
graph without loops and multi-edges, where $V=V(G)$ is the
vertex-set and $E=E(G)$ is the edge-set, where $E\ne\emptyset$.

For $x\in V$, the {\it open neighborhood}, the {\it closed
neighborhood} and the {\it degree} of $x$ are denoted by
$N_{G}(x)=\{y\in V : xy\in E\}$, $N_{G}[x]=N_{G}(x)\cup \{x\}$ and
$deg_G(x)=|N_G(x)|$, respectively. $\delta(G)=\min\{deg_G(x): x\in V\}$
and $\Delta(G)=\max\{deg_G(x): x\in V\}$ are the minimum degree and
the maximum degree of $G$, respectively. For any $X\subseteq V$, let $N_G[X]=\cup_{x\in X}N_G[x]$.

For a subset $D\subseteq V$,  let $\overline D=V\setminus D$. The
notation $G^c$ denotes the complement of $G$, that is , $G^c$ is the
graph with vertex-set $V(G)$ and edge-set $\{xy:\ xy\notin E(G)\
{\rm for\ any}\ x,y\in V(G) \}$. For $B\subseteq E(G^c)$, we use
$G+B$ to denote the graph with vertex-set $V$ and edge-set $E\cup
B$. For convenience, we denote $G+\{xy\}$ by $G+xy$ for an $xy\in
E(G^c)$.

A nonempty subset $D\subseteq V$ is called a {\it dominating set} of
$G$ if $|N_G(x)\cap D|\geq 1$ for each $x\in \overline D$. The {\it
domination number} $\gamma (G)$ of $G$ is the minimum cardinality of
all dominating sets in $G$. The domination is a classical concept in
graph theory. The early literature on the domination with related
topics is, in detail, surveyed in the two books by Haynes,
Hedetniemi, and Slater~\cite{hs981,hs982}.

In 1985, Fink and Jacobson \cite{fj85} introduced the concept of a
generalization domination in a graph. Let $p$ be a positive integer.
A subset $D\subseteq V$ is a {\it $p$-dominating set} of $G$ if
$|N_G(x)\cap D|\geq p$ for each $x\in \overline D$. The {\it
$p$-domination number} $\gamma_p(G)$ is the minimum cardinality of
all $p$-dominating sets in $G$. A $p$-dominating set with
cardinality $\gamma_p(G)$ is called a $\gamma_p$-set of $G$. For $S,
T\subseteq V$, the set $S$ can $p$-dominate $T$ in $G$ if
$|N_G(x)\cap S|\geq p$ for every $x\in T\setminus S$. Clearly, the
$1$-dominating set is the classical dominating set, and so
 $\gamma_1(G)=\gamma(G)$. The $p$-domination is
investigated by many authors (see, for example,
\cite{bcf05,bcv06,cfhv11,cr90,f85}). Very recently, Chellali {\it et
al.}\cite{cfhv11} have given an excellent survey on this topics. The
following are two simple observations.

\begin{obs}\label{obs1.1}
If $G$ is a graph with $|V(G)|\geq p$, then $\g_p(G)\geq p$.
\end{obs}

\begin{obs}\label{obs1.2}
Every $p$-dominating set of a graph contains all vertices of degree at most $p-1$.
\end{obs}

Clearly, addition of some extra edges to a graph could result in
decrease of its domination number. In 1990, Kok and
Mynhardt~\cite{km90} first investigated this problem and proposed
the concept of the reinforcement number. The {\it reinforcement
number} $r(G)$ of a graph $G$ is defined as the smallest number of
edges whose addition to $G$ results in a graph $G'$ with
$\g(G')<\g(G)$. By convention $r(G)=0$ if $\g(G)=1$.

The reinforcement number has received much research attention (see,
for example, \cite{bgh08,dhtv98,hwx09}), and its many variations
have also been well described and studied in graph theory, including
total reinforcement \cite{hrr11,ses07}, independence reinforcement
\cite{zls03}, fractional reinforcement \cite{csm03,dl97} and so on.
In particular, Blair {\it et al.}~\cite{bgh08},
Hu and Xu~\cite{hx10}, independently, showed that the problem
determining $r(G)$ for a general graph $G$ is NP-hard.

Motivated by the work of Kok and
Mynhardt~\cite{km90}, in this paper, we introduce the $p$-reinforcement number, which
is a natural extension of the reinforcement number.
The {\it $p$-reinforcement number} $r_p(G)$ of a graph $G$ is
the smallest number of edges of $G^c$
that have to be added to $G$ in order to
reduce $\gamma_p(G)$, that is
 $$
 r_p(G)=\min\{|B|: B\subseteq E(G^c)\  {\rm with}\ \g_p(G+B)< \g_p(G)\}.
 $$

It is clear that $r_1(G)=r(G)$. By Observation \ref{obs1.1}, we can also
make a convention, $r_p(G)=0$ if $\g_p(G)\leq p$. Thus $r_p(G)$ is well-defined for any graph $G$ and integer $p\geq 1$.
In this paper, we always assume $\g_p(G)> p$ when we consider the $p$-reinforcement number for a graph $G$.

The rest of this paper is organized as follows. In Section 2 we
present an equivalent parameter for calculating the
$p$-reinforcement number of a graph. As its applications, we
determine the values of the $p$-reinforcement numbers for special
classes of graphs such as paths, cycles and complete $t$-partite
graphs in Sections 3, and show that the decision problem on
$p$-reinforcement is NP-hard for a general graph and a fixed integer
$p\geq 2$ in Section 4. Finally, we establish some upper bounds for
the $p$-reinforcement number of a graph $G$ by terms of other
parameters of $G$ in Section 5.

\section{Preliminary}

Let $G$ be a graph with $\g(G)> 1$ and $B\subseteq E(G^c)$
with $|B|=r(G)$ such that $\g(G+B)<\g(G)$. Let $X$ be a $\g$-set
of $G+B$. Then $|B|\geq |V(G)\setminus N_G[X]|$. On the other hand,
given any set $X\subseteq V(G)$, we can always choose a subset
$B\subseteq E(G^c)$ with $|B|=|V(G)\setminus N_G[X]|$ such that $X$
dominates $G+B$. It is a simple observation that, to calculate $r(G)$,
Kok and Mynhardt~\cite{km90} proposed the following parameter
 \begin{equation}\label{e2.1}
 \eta(G)=\min\{|V(G)\setminus N_G[X]|:\ X\subseteq V(G), |X|<\g(G)\},
 \end{equation}
and showed $r(G)=\eta(G)$. We can refine this technique to deal with
the $p$-reinforcement number $r_p(G)$.

Let $G$ be a graph with $\g_p(G)> p$. For any $X\subseteq
V(G)$, let
 \begin{equation}\label{e2.2}
 X^*=\{x\in \overline X: |N_G(x)\cap X|< p\}.
 \end{equation}
Let $B\subseteq E(G^c)$ with $|B|=r_p(G)$ such that
$\g_p(G+B)<\g_p(G)$, and let $X$ be a $\g_p$-set of $G+B$. Then
 $$
 |B|\geq \sum_{x\in X^*}(p-|N_G(x)\cap X|).
 $$
On the other hand, given any set $X\subseteq V(G)$ with $|X|\geq p$, we can always
choose a subset $B\subseteq E(G^c)$ with
 $$
 |B|=\sum_{x\in X^*}(p-|N_G(x)\cap X|)
 $$
such that $X$ can $p$-dominate $G+B$. Motivated by this observation, we
introduce the following notations.
For a subset $X\subseteq V(G)$,
 \begin{eqnarray}
 \eta_p(x,X,G)&=&  \left\{  \begin{array}{ll}
                                     p-|N_G(x)\cap X| & \mbox{if }x\in X^*\\
                                     0                & \mbox{otherwise}
                             \end{array}
                   \right. \mbox{ for $x\in V(G)$,}           \label{e2.3}\\
 \eta_p(S, X, G)&=&\sum_{x\in S}\eta_p(x,X,G) \mbox{\ \ for $S\subseteq V(G)$, and } \label{e2.4} \\
 \eta_p(G)&=&\min\{\eta_p(V(G),X,G) : |X|<\g_p(G)\} \label{e2.5}.
 \end{eqnarray}

A subset $X\subseteq V(G)$ is called an \emph{$\eta_p$-set} of $G$ if
$\eta_p(G)=\eta_p(V(G),X,G)$. Clearly, for any two subsets $S', S\subseteq V(G)$
and two subsets $X', X\subseteq V(G)$,
 $$
 \begin{array}{ll}
 \eta_p(S', X, G)\leq \eta_p(S, X, G) & {\rm if}\ S'\subseteq S,\\
 \eta_p(S, X, G)\leq \eta_p(S, X', G) & {\rm if}\ |X'|\leq |X|.
 \end{array}
 $$
Thus, we have the following simple observation.

\begin{obs}\label{obs2.1}
If $X$ is an $\eta_p$-set of a graph $G$, then $|X|=\g_p(G)-1$.
\end{obs}

The following result shows that computing $r_p(G)$ can be referred to computing $\eta_p(G)$
for a graph $G$ with $\g_p(G)\geq p+1$.

\begin{thm}\label{thm2.2}
For any graph $G$ and positive integer $p$, $r_p(G)=\eta_p(G)$ if
$\g_p(G)> p$.
\end{thm}

\begin{pf}
Let $X$ be an $\eta_p$-set of $G$. Then $|X|=\g_p(G)-1$ by
Observation~\ref{obs2.1}. Let $Y=\{y\in V(G): \eta_p(y,X,G)>0\}$.
Then $Y=X^*$ is contained in $\overline X$, where $X^*$ is defined
in (\ref{e2.2}). Thus, $\eta_p(G)=\eta_p(X^*,X,G)$. We construct a
new graph $G'$ from $G$, for each $y\in X^*$, by adding
$\eta_p(y,X,G)$ edges of $G^c$ to $G$ joining $y$ to $\eta_p(y,X,G)$
vertices in $X$. Clearly, $X$ is a $p$-dominating set of $G'$, that
is, $\g_p(G')\leq |X|$. Let $B=E(G')-E(G)$. Then
 $$
 \g_p(G)=|X|+1>|X|\geq \g_p(G')=\g_p(G+B),
 $$
which implies $r_p(G)\leq |B|$. It follows that
 \begin{equation}\label{e2.6}
 r_p(G)\leq |B|=\sum_{y\in X^*}\eta_p(y,X,G)=\eta_p(X^*,X,G)=\eta_p(G).
 \end{equation}

On the other hand, let $B$ be a subset of $E(G^c)$ such that
$|B|=r_p(G)$ and $\g_p(G+B)=\g_p(G)-1$. Let $G'=G+B$ and $X'$ be a
$\g_p$-set of $G'$. For every $xy\in B$,  $X'$ cannot $p$-dominate
the graph $G'-xy$ by the minimality of $B$. This fact means that
only one of $x$ and $y$ is in $X'$. Without loss of generality,
assume $y\in \overline {X'}$. Since $X'$ cannot $p$-dominate $y$ in
$G'-xy$ and so in $G$, $|N_G(y)\cap X'|<p$. Let $Z$ be all
end-vertices of edges in $B$ and $Y=\overline{X'}\cap Z$.
Since $X'$ is a $\g_p$-set of
$G'$, $|N_{G'}(u)\cap X'|\geq p$ for any $u\in \overline{X'}$. In
other words, any $u\in \overline{X'}$ with $|N_{G}(u)\cap X'|<p$
must be in $Y$. It follows that
 \begin{equation}\label{e2.7}
 \sum_{u\in \overline{X'}}\eta_p(u,X',G)=\sum_{y\in Y}(p-|N_G(y)\cap X'|)=|B|.
 \end{equation}
By (\ref{e2.7}), we immediately have that
 $$
\eta_p(G)\leq \eta_p(V(G),X',G)=\sum_{u\in \overline{X'}}\eta_p(u,X',G)=|B|=r_p(G).
 $$
Combining this with (\ref{e2.6}), we obtain $r_p(G)=\eta_p(G)$, and
so the theorem follows.
\end{pf}

Note that when $p=1$, $X^*$ defined in (\ref{e2.2}) is
$V(G)\setminus N_G[X]$. This fact means that $\eta(G)$ defined in
(\ref{e2.1}) is a special case of $p=1$ in (\ref{e2.5}), that is,
$\eta_1(G)=\eta(G)$. Thus, the following corollary holds
immediately.

\begin{cor}\label{cor2.1}\textnormal{(Kok and Mynhardt~\cite{km90})}
$r(G)=\eta(G)$  if $\g(G)>1$.
\end{cor}

Using Observation \ref{obs1.2} and Theorem \ref{thm2.2}, the
following corollary is obvious.

\begin{cor}\label{cor2.2}
Let $p\geq 1$ be an integer and $G$ be a graph with $\g_p(G)> p$. If $\Delta(G)< p$, then
$$r_p(G)=p-\Delta(G).$$
\end{cor}

\section{Some Exact Values}

In this section we will use Theorem \ref{thm2.2} to calculate the
$p$-reinforcement numbers for some classes of graphs.

We first determine the $p$-reinforcement numbers for paths and
cycles. Let $P_n$ and $C_n$ denote, respectively, a path and a cycle
with $n$ vertices. When $p=1$, Kok and Mynhardt \cite{km90} proved
that $r(P_n)=r(C_n)=i$ if $n=3k+i \geq 4$, where $i\in \{1,2,3\}$.
We will give the exact values of $r_p(P_n)$ and $r_p(C_n)$ for
$p\geq 2$. The following observation is simple but useful.

\begin{obs}\label{obs3.1}
For integer $p\geq 2$,
$$
\g_p(P_n)=\left\{
              \begin{array}{rl}
              \lfloor\frac{n}{2}\rfloor+1 & \mbox{\ \ if\ \ $p=2$}\\
              n                       & \mbox{\ \ if\ \ $p\geq 3$}
              \end{array}
           \right.
\mbox{ and\ \ }
\g_p(C_n)=\left\{
              \begin{array}{rl}
              \lceil\frac{n}{2}\rceil & \mbox{\ \ if\ \ $p=2$}\\
              n                       & \mbox{\ \ if\ \ $p\geq 3$.}
              \end{array}
           \right.
$$
\end{obs}

\begin{thm}\label{thm3.2}
Let $p\geq 2$  be an integer. If $\g_p(P_n)>p$ then
$$
r_p(P_n)=\left\{
              \begin{array}{ll}
              2 & \mbox{\ \ if\ \ $p=2$ and $n$ is odd}\\
              1 & \mbox{\ \ if\ \ $p=2$ and $n$ is even}\\
              p-2                      & \mbox{\ \ if\ \ $p\geq 3$}.
              \end{array}
           \right.
$$
\end{thm}

\begin{pf}
Let $P_n=x_1x_2\cdots x_n$ and $X$ be an $\eta_p$-set of $P_n$. By
Theorem \ref{thm2.2} and $\g_p(P_n)> p$,
$r_p(P_n)=\eta_p(P_n)=\eta_p(V(P_n),X,P_n)\geq 1$. For $p\geq 3$, it
is easy to see that $r_p(P_n)=p-2$ by Corollary \ref{cor2.2}. Assume
that $p=2$ below.

If $n$ is even, then by Observation \ref{obs3.1},
$\g_2(P_n)-\g_2(C_n)=1$, which implies that $r_2(P_n)\leq 1$.
Furthermore, $r_2(P_n)=1$.

If $n$ is odd, then $\g_2(P_n)=\frac{n+1}{2}$ by Observation
\ref{obs3.1}, and so $n\geq 5$ since $\g_2(P_n)>2$. Let
 $$
 X'=\bigcup_{i=1}^{\frac{n-1}{2}}\{x_{2i}\}.
 $$
Clearly, $|X'|=\frac{n-1}{2}=\g_2(P_n)-1$. So
\begin{eqnarray*}
 \eta_2(V(P_n),X,P_n)\leq \eta_2(V(P_n),X',P_n)
                                 = \eta_2(x_1, X',P_n)+\eta_2(x_n,X',P_n)
                                 = 2.
\end{eqnarray*}

Suppose that $\eta_2(V(P_n),X,P_n)=1$. Then $X$ can $2$-dominate
either $V(P_n)\setminus \{x_1\}$ or $V(P_n)\setminus \{x_{n}\}$. In
both cases, we have
 $$
 |X|\geq \g_2(P_{n-1})=\left\lfloor\frac{n-1}{2}\right\rfloor+1=\frac{n-1}{2}+1,
 $$
which contradicts with $|X|=\frac{n-1}{2}$. Hence
$r_2(P_n)=\eta_2(V(P_n),X,P_n)=2$.
\end{pf}

\begin{thm}\label{3.3}
Let $p\geq 2$ be an integer. If $\g_p(C_n)>p$ then
$$
r_p(C_n)=\left\{
              \begin{array}{ll}
              2 & \mbox{\ \ if\ \ $p=2$ and $n$ is odd}\\
              4 & \mbox{\ \ if\ \ $p=2$ and $n$ is even}\\
              p-2                      & \mbox{\ \ if\ \ $p\geq 3$}.
              \end{array}
           \right.
$$
\end{thm}

\begin{pf}
Let $C_n=x_1x_2\cdots x_nx_1$. If $p\geq 3$ then the result holds
obviously by Corollary \ref{cor2.2}. In the following, we only need
to calculate the values of $r_p(C_n)$ for $p=2$. Let $X$ be an
$\eta_2$-set of $C_n$. Then
$r_2(C_n)=\eta_2(C_n)=\eta_2(V(C_n),X,C_n)$ by Theorem \ref{thm2.2}.
Note that $n\geq 5$  since $\g_2(C_n)=\lceil\frac{n}{2}\rceil>2$.

If $n$ is odd, then let
 $$
 X'=\bigcup_{i=1}^{\frac{n-1}{2}}\{x_{2i-1}\}.
 $$
Clearly, $|X'|=\frac{n-1}{2}=\g_2(C_n)-1$ by Observation
\ref{obs3.1}, and
$\eta_2(V(C_n),X',C_n)=\eta_2(x_{n-1},X',C_n)+\eta_2(x_n,X',C_n)=2$.
So
 $$
 r_2(C_n)=\eta_2(V(C_n),X,C_n)\leq \eta_2(V(C_n),X',C_n)=2.
 $$
Since $X$ is not a $2$-dominating set of $C_n$,  there must be two
adjacent vertices, denoted by $x_i$ and $x_{i+1}$, of $C_n$ not in
$X$. This fact means that $\eta_2(x_i,X,C_n)\geq 1$ and
$\eta_2(x_{i+1},X,C_n)\geq 1$. So
$$r_2(C_n)=\eta_2(V(C_n),X,P_n)\geq \eta_2(x_i,X,C_n)+\eta_2(x_{i+1},X,C_n)\geq 2.$$
Hence $r_2(C_n)=2$.

If $n$ is even, then $n\geq 6$. Deleting $X$ and all vertices
$2$-dominated by $X$ from $C_n$, we can obtain a result graph,
denoted by $H$, each of whose components is a path with length at
least 2. Denote all components of $H$ by $H_1,\cdots, H_h$, where
$h\geq 1$. In the case that $h=1$ and the length of $H_1$ is equal
to one,  $X$ can $2$-dominate a subgraph of $C_n$ that is isomorphic
to $P_{n-2}$. By Observation \ref{obs3.1},
 $$
 |X|\geq \g_2(P_{n-2})=\lfloor\frac{n-2}{2}\rfloor+1=\frac{n}{2},
 $$
which contradicts that
$|X|=\g_2(C_n)-1=\lceil\frac{n}{2}\rceil-1=\frac{n}{2}-1$. In other
cases, we can find that
 $$
 r_2(C_n)=\eta_2(V(C_n),X,C_n)\geq 4.
 $$
Let
 $$
 X''=\bigcup_{i=1}^{\frac{n}{2}-1}\{x_{2i-1}\}.
 $$
It is easy to check that $|X''|=\frac{n}{2}-1=\g_2(C_n)-1$ and
$\eta_2(V(C_n),X'',C_n)=4$. So
$$r_2(C_n)=\eta_2(V(C_n),X,C_n)\leq\eta_2(V(C_n),X'',C_n)=4.$$
Hence $r_2(C_n)=4$ and so the theorem is true.
\end{pf}

Next we consider the $p$-reinforcement number for a complete
$t$-partite graph $K_{n_1,\cdots,n_t}$. To state our results, we
need some symbols. For any subset $X=\{n_{i_1,\ \cdots,\ n_{i_r}}\}$
of $\{n_1, \cdots, n_t\}$, define
 $$
 |X|=r \mbox{\ \ and\ \ } f(X)=\sum_{j=1}^rn_{i_j}.
 $$
For convenience,  let $|X|=0$ and $f(X)=0$ if $X=\emptyset$. let
$$
\mathscr{X}=\{X : X \mbox{ is a subset of } \{n_1, \cdots, n_t\}
\mbox{ with $f(X)\geq \g_p(G)$}\}
$$
and, for every $X\in \mathscr{X}$, define
$$
f^*(X)=\max\{f(Y) : Y \mbox{ is a subset of $X$ with $|Y|=|X|-1$ and
$f(Y)< p$}\}.
$$

\begin{thm}\label{thm3.4}
For any integer $p\geq 1$ and a complete $t$-partite graph
$G=K_{n_1, \cdots, n_t}$  with $t\geq 2$ and $\g_p(G)> p$,
$$
r_p(G)=\min\{(p-f^*(X))(f(X)-\g_p(G)+1) : X\in \mathscr{X}\}.
$$
\end{thm}

\begin{pf}
Let $N=\{n_1,\cdots,n_t\}$ and $V(G)=V_1\cup \cdots \cup V_t$ be the
vertex-set of $G$ such that $|V_i|=n_i$ for each $i=1,\cdots,t$. Let
 $$m=\min\{(p-f^*(X))(f(X)-\g_p(G)+1) : X\in \mathscr{X}\}.$$

We first prove that $r_p(G)\leq m$. Let $X\subseteq \mathscr{X}$
(without loss of generality, assume $X=\{n_1,\cdots,n_k,n_{k+1}\}$
for some $0\leq k\leq t-1$) such that
\begin{eqnarray*}
f^*(X)=n_1+\cdots+n_k \mbox{ and } (p-f^*(X))(f(X)-\g_p(G)+1)=m.
\end{eqnarray*}
By $X\subseteq \mathscr{X}$, we know that
$n_{k+1}=f(X)-f^*(X)\geq\g_p(G)-f^*(X)$.  So we can pick a
vertex-subset $V_{k+1}'$ from $V_{k+1}$ such that
$|V_{k+1}'|=\g_p(G)-f^*(X)-1$. Let
 $$
 D=V_1\cup \cdots \cup V_k\cup V_{k+1}'.
 $$
Clearly,  $|D|=\g_p(G)-1.$ Since $\g_p(G)> p$,  $|D|\geq p$ and so
$D$ can $p$-dominate  $\cup_{i=k+2}^tV_i$. Hence by the definition
of $\eta_p(V(G),D,G)$,
\begin{eqnarray*}
\eta_p(V(G),D,G)&=&\eta_p(V(G)\setminus D, D,G)\\
                &=&\sum_{v\in V_{k+1}\setminus V_{k+1}'}\eta_p(v,D,G)+\sum_{i=k+2}^t\eta_p(V_i,D,G)\\
                &=& |V_{k+1}\setminus V_{k+1}'|(p-f^*(X))+0\\
                &=&(p-f^*(X))[n_{k+1}-(\g_p(G)-f^*(X)-1)]\\
                &=&(p-f^*(X))(f(X)-\g_p(G)+1)\\
                &=& m.
\end{eqnarray*}
By Theorem \ref{thm2.2}, we have
$r_p(G)=\eta_p(G)\leq \eta_p(V(G),D,G)=m.$

On the other hand, we will show that $r_p(G)\geq m$. For any subset
$M$ of $N$, we use $I(M)$ to denote the subindex-sets of all
elements in $M$, that is,
$$I(M)=\{i : n_i\in M\}.$$

 Let $S$ be an $\eta_p$-set of $G$ and
let
\begin{eqnarray*}
&& Y=\{n_i : \hspace{0.7cm}|V_i\cap S|=|V_i| \mbox{ for $1 \leq i\leq t$}\}, \mbox{ and}\\
&& A=\{n_i : 0< |V_i\cap S|<|V_i| \mbox{ for $1 \leq i\leq t$}\}.
\end{eqnarray*}
Thus
\begin{equation}\label{e5}
f(Y\cup A)=f(Y)+f(A)=\sum_{i\in I(Y)}|V_i|+\sum_{i\in I(A)}|V_i|\\
                    \geq |S|=\g_p(G)-1
\end{equation}
by Observation \ref{obs2.1}. Since $\cup_{i\in I(Y)}V_i\ (\subseteq
S)$ cannot $p$-dominate $G$,
\begin{equation}\label{e6}
f(Y)=\sum_{i\in I(Y)}n_i=|\cup_{i\in I(Y)}V_i|< p.
\end{equation}
Hence,  by (\ref{e5}) and $\g_p(G)> p$,
$$f(A)\geq \g_p(G)-1-f(Y)> \g_p(G)-p-1\geq 0,$$
 which implies that $|A|\geq 1$.

\noindent \textbf{Claim.}  $|A|=1$.

\noindent \textbf{Proof of Claim.} Suppose that $|A|\geq 2$. Then we
can choose $i$ and $j$ from $I(A)$ such that $i\neq j$. By the
definition of $A$, we have $0< |V_i\cap S|< |V_i|$ and $0< |V_j\cap
S|< |V_j|$. Therefore, we can pick two vertices $x$ and $y$ from
$V_i\cap S$ and $V_j\setminus S$, respectively. Let
 $$
 S'=(S\setminus \{x\})\cup \{y\}.
 $$
Obviously, $|S'|=|S|=\g_p(G)-1$, $|V_i \cap S'|=|V_i\cap S|-1$ and
$|V_j \cap S'|=|V_j\cap S|+1$.

Note that $G$ is a complete $t$-partite graph. For any $v\in V(G)$,
we can easily find the value of $\eta_p(v,S',G)-\eta_p(v,S,G)$ by
the definitions of $\eta_p(v,S',G)$ and $\eta_p(v,S,G)$ as follows:
\begin{equation*}
\eta_p(v,S',G)-\eta_p(v,S,G)=
     \left\{
       \begin{array}{rl}
       (p-|S|+|V_i\cap S|-1)-0   & \ \ \ \mbox{if }v=x\\
       -1                        & \ \ \ \mbox{if }v\in V_i\setminus S\\
       0-(p-|S|+|V_j\cap S|)     & \ \ \ \mbox{if }v=y\\
       1                         & \ \ \ \mbox{if }v\in (V_j\setminus S)\setminus \{y\}\\
       0                         & \ \ \ \mbox{otherwise.}
       \end{array}
     \right.
\end{equation*}
Since $S$ is an $\eta_p$-set of $G$ and $|S'|=|S|$, we have
\begin{eqnarray*}
0&\leq& \eta_p(V(G),S',G)-\eta_p(V(G),S,G)\\
  &=& \sum_{v\in V(G)}(\eta_p(v,S',G)-\eta_p(v,S,G))\\
  &=& (p-|S|+|V_i\cap S|-1)-|V_i\setminus S|-(p-|S|+|V_j\cap S|)+|(V_j\setminus S)\setminus \{y\}|\\
  &=& (|V_i\cap S|-|V_i\setminus S|)-(|V_j\cap S|-|V_j\setminus S|)-2.
\end{eqnarray*}
This means that
$$(|V_i\cap S|-|V_i\setminus S|)\geq(|V_j\cap S|-|V_j\setminus S|)+2.$$
However, by the symmetry of $V_i$ and $V_j$, we can also obtain
$$(|V_j\cap S|-|V_j\setminus S|)\geq(|V_i\cap S|-|V_i\setminus S|)+2$$
 by applying the similar discussion. This is a contradiction, and so the claim holds. $\Box$

By \textbf{Claim}, we can assume that $I(A)=\{h\}$. From the
definitions of $Y$ and $A$, we have $|Y\cup A|=|Y|+1$ and
  \begin{equation*}
  f(Y\cup A)=\sum_{i\in I(Y)}|V_i|+|V_h|\geq \sum_{i\in I(Y)}|V_i|+(|V_h\cap S|+1)= |S|+1=\g_p(G).
  \end{equation*}
It follows that $Y \cup A \in \mathscr{X}$. Thus, by (\ref{e6}) and
the definition of $f^*(Y\cup A)$,  we have $f(Y)\leq f^*(Y\cup A)$.
Since $\g_p(G)>p$, $|S|=\g_p(G)-1\geq p$, and so $S$ $p$-dominates
$V(G)\setminus (\cup_{i\in I(Y\cup A)}V_i)$.
 Therefore, by Theorem \ref{thm2.2},
\begin{eqnarray*}
r_p(G)=\eta_p(G)=\eta_p(V(G),S,G)&=&\eta_p(V(G)\setminus S,S,G)\\
                &=&\sum_{v\in V_h\setminus S}\eta_p(v,S,G)\\
                &=& (p-f(Y))|V_h\setminus S|\\
                &=& (p-f(Y))[|V_h|-(|S|-f(Y))]\\
                &=& (p-f(Y))(f(Y\cup A)-\g_p(G)+1)\\
                &\geq&(p-f^*(Y\cup A))(f(Y\cup A)-\g_p(G)+1)\\
                &\geq& m.
\end{eqnarray*}
This completes the proof of the theorem.
\end{pf}

\vskip6pt

For example, let $G=K_{2,2,10,17}$ and $p=11$. Then $\g_{11}(G)=12$,
and so
$$
\mathscr{X}=\{\{17\},\{2,10\},\{2,17\},\{10,17\},\{2,2,10\},\{2,2,17\},\{2,10,17\},
                                \{2,2,10,17\}\}.
$$
By Theorem~\ref{thm3.4}, for any $X\in \mathscr{X}$, we have that
$$
f^*(X)=\left\{
              \begin{array}{rl}
              0  &   \mbox{ if $X=\{17\},\{2,10,17\}$ or $\{2,2,10,17\}$};\\
              2  &    \mbox{ if $X=\{2,17\}$};\\
              4  &     \mbox{ if $X=\{2,2,10\}$ or $\{2,2,17\}$};\\
              10 &      \mbox{ if $X=\{2,10\}$ or $\{10,17\}$}.
              \end{array}
      \right.
$$
Hence
\begin{eqnarray*}
r_{11}(G)&=&\min\{(11-f^*(X))(f(X)-\g_{11}(G)+1): X\in \mathscr{X}\}\\
         &=& \min\{(11-f^*(X))(f(X)-11): X\in \mathscr{X}\}\\
         &=& (11-f^*(\{2,10\}))(f(\{2,10\})-11)\\
         &=& 1.
\end{eqnarray*}

\section{Complexity}

Blair et al. \cite{bgh08}, Hu and Xu \cite{hx10}, independently,
showed that the $1$-reinforcement problem is NP-hard. Thus, for any
positive integer $p$, the $p$-reinforcement problem is also NP-hard
since the $1$-reinforcement is a sub-problem of the
$p$-reinforcement problem.

For each fixed $p$, $p$-dominating set is polynomial-time computable
(see Downey and Fellows~\cite{df95,df97} for definitions and
discussion). However, the $p$-reinforcement number problem is hard
even for specific values of the parameters. In this section, we will
consider the following decision problem.

\textbf{$p$-Reinforcement}

\emph{Instance}: A graph $G$, $p\ (\geq 2)$ is a fixed integer.

\emph{Question}: Is $r_p(G)\leq 1$?

We will prove that \textbf{$p$-Reinforcement} ($p\geq 2$) is also
NP-hard by describing a polynomial transformation from the following
NP-hard problem (see \cite{gj79}).

\noindent\begin{eqnarray*}
&&\mbox{\textbf{3-Satisfiability (3SAT)}}\\
&&\mbox{\emph{Instance}: A set $U =\{u_1,\ldots, u_n\}$ of variables and a collection $\mathscr{C} = \{C_1,\ldots,C_m\}$}\\
&&\mbox{\hspace{1.75cm} of clauses over $U$ such that $|C_i | = 3$ for $i = 1, 2,\ldots,m$. } \\
&&\mbox{\hspace{1.75cm} Furthermore, every literal is used in at least one clause.}\\
&&\mbox{\emph{Question}: Is there a satisfying truth assignment for $C$?}
\end{eqnarray*}

\begin{thm}\label{t2}
 For a fixed integer $p\geq 2$, \textbf{$p$-Reinforcement} is NP-hard.
\end{thm}

\begin{pf}
Let $U =\{u_1,\ldots, u_n\}$ and $\mathscr{C} =\{C_1,\ldots,C_m\}$
be an arbitrary instance $I$ of \textbf{3SAT}.
We will show the NP-hardness of \textbf{$p$-Reinforcement} by reducing \textbf{3SAT} to it in polynomial time.
To this aim,  we construct a graph $G$ as follows:
\begin{description}
  \item[\quad a.] For each variable $u_i\in U$, associate a graph $H_i$, where $H_i$ can be obtained
from a complete graph $K_{2p+2}$ with vertex-set $\{u_i,\overline  u_i\}\cup(\cup_{j=1}^p\{v_{i_j},\overline v_{i_j}\})$ by
deleting the edge-subset $\cup_{j=1}^{p-1}\{u_i \overline v_{i_j},\overline u_i v_{i_j}\}$;
  \item[\quad b.] For each clause $C_j\in \mathscr{C}$, create a single vertex $c_j$ and join $c_j$ to the vertex $u_i$ (resp. $\overline{u}_i$) in $H_i$ if and only if the literal $u_i$ (resp. $\overline{u}_i$) appears in clause $C_j$ for any $i\in \{1,\ldots,n\}$;
  \item[\quad c.] Add a complete graph $T\ (\cong K_p)$ and join all of its vertices to each $c_j$.
\end{description}

For convenience, let $X_i=\cup_{j=1}^p\{v_{i_j}\}$ and $\overline
X_i=\cup_{j=1}^p\{\overline v_{i_j}\}$. Then $V(H_i)=\{u_i,\overline
u_i\}\cup X_i\cup \overline X_i$. Use $H_0$ to denote the induced
subgraph by $\{c_1,\cdots,c_m\}\cup V(T)$.

It is clear that the construction of $G$ can be accomplished in
polynomial time. To complete the proof of the theorem, we only need
to prove that $\mathscr{C}$ is satisfiable if and only if
$r_p(G)=1$. We first prove the following two claims.

\noindent\textbf{Claim 1.} {\it Let $D$ be a $\g_p$-set of $G$. Then
$|D|=p(n+1)$, moreover, $|V(H_i)\cap D|=p$ and $|\{u_i,\overline
u_i\}\cap D|\leq 1$ for each $i\in \{1,2,\ldots,n\}$.}

\noindent {\bf Proof of Claim 1.} Suppose there is some $i\in
\{1,2,\cdots,n\}$ such that $|V(H_i)\cap D|<p$. Then there must be a
vertex, say $x$, of $V(H_i)\setminus D$ such that $N_G(x)\subseteq
V(H_i)$. And so $|N_G(x)\cap D|\leq |V(H_i)\cap D|< p$, which
contradicts that $D$ is a $\g_p$-set of $G$. Thus $|V(H_i)\cap
D|\geq p$ for each $i\in \{0,1,\cdots,n\}$, and so
  \begin{equation}\label{eq4.1}
  \g_p(G)=|D|=\sum_{i=0}^n|V(H_i)\cap D|\geq p(n+1).
  \end{equation}
On the other hand, let
 $$
 D'=\bigcup_{i=1}^n[(X_i- \{v_{i_p}\})\cup \{\overline u_i\}]\cup V(T).
 $$
Clearly, $|D'|=p(n+1)$ and $D'$ is a $p$-dominating set of $G$. Hence by (\ref{eq4.1}),
 $$
 p(n+1)\leq \sum_{i=0}^n|V(H_i)\cap D|= \g_p(G)\leq |D'|= p(n+1),
 $$
which implies that $\g_p(G)=p(n+1)$ and $|V(H_i)\cap D|=p$ for each $0\leq i\leq n$. Furthermore,
 if $|\{u_i,\overline u_i\}\cap D|=2$ then $|(X_i\cup \overline X_i)\cap D|=p-2$. So we can choose a vertex from $X_i\cup \overline X_i$ that is not $p$-dominated by $D$. This is impossible since $D$ is a $\g_p$-set of $G$, and so $|\{u_i,\overline u_i\}\cap D|\leq 1$.  The claim holds. $\square$

\noindent\textbf{Claim 2.}
{\it If there is an edge $e=xy\in G^c$ such that $\g_p(G+e)<\g_p(G)$, then any $\g_p$-set $D_e$ of $G+e$ satisfies the following properties.
\begin{description}
  \item[\ \, $(i)$] $|V(H_i) \cap D_e|=p$ and $|\{u_i,\overline u_i\}\cap D_e|\leq 1$ for each $i\in \{1,\cdots,n\}$;
  \item[\, $(ii)$] $\{c_1,\cdots,c_m\}\cap D_e=\emptyset$, and so $|V(T)\cap D_e|=p-1$;
  \item[\,$(iii)$] One of $x$ and $y$ belongs to $V(T)\setminus D_e$ and the other belongs to $H\cap D_e$, where $H=\cup_{i=1}^nV(H_i)$.
\end{description}
}
\textbf{Proof of Claim 2.} Because $D_e$ is a $\g_p$-set of $G+e$ and $\g_p(G+e)<\g_p(G)$,  one of $x$ and $y$ is not in $D_e$ but the other is in $D_e$. Without loss of generality, say $x\notin D_e$ and $y\in D_e$. It is clear that $|N_G(x)\cap D_e|=p-1$. Since vertex $x$ is the unique vertex not be $p$-dominated by $D_e$, we have
\begin{equation}\label{eq4.2}
\eta_p(V(G),D_e,G)=\eta_p(x,D_e,G)=p-(p-1)=1.
\end{equation}
 Let
 $$
 D=D_e\cup \{x\}.
 $$
Then $D$ is a $p$-dominating set of $G$ and
$|D|=|D_e|+1=\g_p(G+e)+1\leq\g_p(G)$. That is, $D$ is a $\g_p$-set
of $G$. By Claim 1,
\begin{equation}\label{eq4.3}
|V(H_i)\cap D|=p \mbox{ for each $i=0,1,\cdots,n$},
\end{equation}
and $|\{u_i,\overline u_i\}\cap D_e|\leq |\{u_i,\overline u_i\}\cap D|\leq 1$ for $1\leq i\leq n$.

Suppose that there exists some $i\in \{1,\cdots,n\}$ such that
$|V(H_i)\cap D_e|\neq p$. Then by (\ref{eq4.3}), $x\in V(H_i)$ and
$|V(H_i)\cap D_e|=p-1$. Thus every vertex in $(X_i\cup \overline
X_i)\setminus (D_e\cup\{x\})$ is dominated by at most $p-1$ vertices
of $D_e$. Hence by $|X_i\cup \overline X_i|=2p$,
 $$
 \eta_p(V(G),D_e,G)\geq\eta_p(X_i\cup \overline X_i,D_e,G)
 \geq |(X_i\cup \overline X_i)\setminus D_e|-1\geq 2p-(p-1)-1>1,
 $$
which contradicts with (\ref{eq4.2}). Hence $(i)$ holds.

Suppose that there is some $j\in \{1,\cdots,m\}$ such that $c_j\in
D_e$. By $(i)$ and (\ref{eq4.3}), $x\in V(H_0)$ and so $|V(H_0)\cap
D_e|=|V(H_0)\cap D|-1=p-1$. Hence $|V(T)\cap D_e|\leq p-2$ by
$V(H_0)=\{c_1,\cdots,c_m\}\cup V(T)$. Since each vertex of $T\
(\cong K_p)$ has exact $p-1$ neighbors in $D_e$,
 $$
 \eta_p(V(G),D_e,G)\geq \eta_p(V(T),D_e,G)=|V(T)\setminus D_e|=p-|V(T)\cap D_e|\geq 2.
 $$
This contradicts with (\ref{eq4.2}). Thus $\{c_1,\cdots,c_m\}\cap
D_e=\emptyset$, and so $|V(T)\cap D_e|=|V(H_0)\cap D_e|=p-1$. Hence
$(ii)$ holds.

By $(ii)$, $T$ has a unique vertex, say $z$, not in $D_e$. From
$|N_G(z) \cap D_e|=|V(H_0)\cap D_e|=p-1$, the vertex $z$ is not
$p$-dominated by $D_e$. However, $x$ is the unique vertex not be
$p$-dominated by $D_e$ in $G$ by (\ref{eq4.2}). Thus $z=x$, and so
$x=z\in V(T)\setminus D_e$. By the construction of $G$ and $xy\in
G^c$, it is clear that $y\in (\cup_{i=1}^nV(H_i))\cap D_e$. Hence
$(iii)$ holds. $\Box$

We now show that $\mathscr{C}$ is satisfiable if and only if
$r_p(G)=1$.

If $\mathscr{C}$ is satisfiable, then $\mathscr{C}$ has a satisfying truth assignment $t: U\rightarrow \{T,F\}$. According to this satisfying assignment, we can choose a subset $S$ from $V(G)$ as follows:
$$S=S_0\cup S_1\cup\cdots\cup S_n,$$
where $S_0$ consists of $p-1$ vertices of $T$ and
$$
S_i=\left\{
           \begin{array}{ll}
           u_i\cup (\overline X_i- \{\overline v_{i_p}\}) & \mbox{ if $t(u_i)=T$}\\
           \overline u_i\cup (X_i- \{v_{i_p}\}) & \mbox{ if $t(u_i)=F$}
           \end{array}
           \mbox{ for each $i\in \{1,\cdots,n\}$.}
     \right.
$$
It can be verified easily that $|S|=p(n+1)-1=\g_p(G)-1$ and
$\cup_{i=1}^nV(H_i)$ can be $p$-dominated by $S$. Since $t$ is a
satisfying true assignment for $\mathscr{C}$, each clause $C_j\in
\mathscr{C}$ contains at least one true literal. That is, the
corresponding vertex $c_j$ has at least one neighbor in $\{u_1,\br
u_1\cdots,u_n,\bar u_n\}\cap S$ by the definitions of $G$ and $S$,
and so every $c_j\in \{c_1,\cdots,c_m\}$ has at least $p$ neighbors
in $S$ since $S_0\subseteq N_G(c_j)$. Note that the unique vertex in
$V(T)\setminus S_0$ has exact $p-1$ neighbors in $S$. By Theorem
\ref{thm2.2} and $|S|=\g_p(G)-1$,
 $$
 r_p(G)=\eta_p(G)\leq \eta_p(V(G),S,G)=\eta_p(V(T)\setminus S_0,S,G)=p-(p-1)=1.
 $$
Furthermore, we have $r_p(G)=1$ since $\g_p(G)>p$ by Claim 1.

Conversely, assume $r_p(G)= 1$.  That is, there exists an edge $e=xy$ in
$G^c$ such that $\g_p(G+e)<\g_p(G)$. Let $D_e$ be a $\g_p$-set of
$G+e$.
Define $t: U\to \{T,F\}$ by
 \begin{equation}\label{eq4.4}
 t(u_i)=\left\{
\begin{array}{ll}
 T \ & \mbox{ if vertex}\ u_i\in D_e \\
 F \ & \mbox{ if vertex}\ u_i\notin D_e
\end{array}
 \right.
 \ \mbox{for }i=1,\cdots,n.
 \end{equation}

We will show that $t$ is a satisfying truth assignment for
$\mathscr{C}$. Let $C_j$ be an arbitrary clause in $\mathscr{C}$.
 By $(ii)$ and $(iii)$ of Claim 2, the corresponding vertex $c_j$ is not in $D_e$ and
$|N_G(c_j)\cap D_e|\geq p$ since $c_j\notin \{x,y\}$. Then there must be some $i\in \{1,\cdots,n\}$ such that
\begin{equation}\label{eq4.5}
|\{u_i,\overline u_i\}\cap N_G(c_j)\cap D_e|=1,
\end{equation}
since $T$ contains exact $p-1$ vertices of $D_e$ by $(i)$ and $(ii)$ of Claim 2.
If $u_i\in N_G(c_j)\cap D_e$, then $u_i\in C_j$ and $t(u_i)=T$ by the construction of $G$ and (\ref{eq4.4}).
If $\overline u_i\in N_G(c_j)\cap D_e$, then the literal $\overline u_i$ belongs to $C_j$ by the construction of $G$. Note that $u_i\notin D_e$ from $\overline u_i\in D_e$ and $(i)$ of Claim 2. This means that $t(u_i)=F$ by (\ref{eq4.4}). Hence $t(\overline u_i)=T$.
The
arbitrariness of $C_j$ with $1\le j\le m$ shows that all the clauses
in $\mathscr{C}$ is satisfied by $t$. That is, $\mathscr{C}$ is
satisfiable.

The theorem follows.
\end{pf}

\section{Upper Bounds}

For a graph $G$ and $p=1$, Kok and Mynhardt \cite{km90} provided an
upper bound for $r(G)$ in terms of the smallest private neighborhood
of a vertex in some $\g$-set of $G$.
Let $X\subseteq V(G)$ and $x\in X$.
The {\it private neighborhood} of $x$ with respect to $X$ is defined
as the set
 \begin{equation}\label{e4.1}
PN(x,X,G)=N_G[x]\setminus N_G[X\setminus\{x\}].
\end{equation}
Set
 $$
 \mu(X,G)=\min\{|PN(x,X,G)|:\ x\in X\}
 $$ and
 \begin{equation}\label{e4.2}
 \mbox{$\mu(G)=\min\{\mu(X,G):\ X$ is a $\g$-set of $G\}$.}
 \end{equation}
Using this parameter, Kok and Mynhardt \cite{km90} showed that
$r(G)\leq \mu(G)$ if $\g(G)\geq 2$ with equality if $\g(G)=1$. We
generalize this result to any positive integer $p$.

In order to state
our results, we need some notations. Let $X\subseteq V(G)$ and $x\in
X$. A vertex $y\in \overline X$ is called a {\it $p$-private
neighbor} of $x$ with respect to $X$ if $xy\in E(G)$ and
$|N_G(y)\cap X|=p$.  The {\it $p$-private
neighborhood} of $x$ with respect to $X$ is defined as
  \begin{equation}\label{e4.3}
 PN_p(x,X,G)=\{y:\ y \mbox{ is a $p$-private neighbor of $x$ with respect to $X$}\}.
 \end{equation}
Let
 \begin{eqnarray}
 \mu_p(x,X,G)&=&|PN_p(x,X,G)|+\max\{0,p-|N_G(x)\cap X|\},\label{e4.4}\\
 \mu_p(X,G)&=&\min\{\mu_p(x,X,G) : x\in X\}, \mbox{\ \ and}\label{e4.5}\\
 \mu_p(G)&=&\min\{\mu_p(X,G):\ X \mbox{ is a $\g_p$-set of $G$}\}\label{e4.6}.
\end{eqnarray}

\begin{thm}\label{thm4.1}
For any graph $G$ and positive integer $p$,
 $$
 r_p(G)\leq \mu_p(G)
 $$
with equality if $r_p(G)=1$.
\end{thm}

\begin{pf}
If $\g_p(G)\leq p$, then $r_p(G)=0\leq \mu_p(G)$ by our convention.
Assume that $\g_p(G)\geq p+1$ below. Let $X$ be a $\g_p$-set of $G$
and $x\in X$ such that
 $$
 \mu_p(G)=\mu_p(X,G)=\mu_p(x,X,G).
 $$
Since
$|X|=\g_p(G)\geq p+1$, we can choose a vertex, say $u_y$, from
$X\setminus N_G(y)$ for each $y\in PN_p(x,X,G)$, and a subset $X'$ with $|X'|=\max\{0,p-|N_G(x)\cap
X|\}$ from $X\setminus N_G[x]$.
Let
 $$
 G'=G+\{yu_y :\ y\in PN_p(x,X,G)\}+\{xv :\ v\in X'\}.
 $$
Obviously, $X\setminus\{x\}$ is a $p$-dominating set of $G'$, which implies that
 $$
 r_p(G)\leq |PN_p(x,X,G)|+|X'|=\mu_p(x,X,G)=\mu_p(G).
 $$

Assume $r_p(G)=1$. Then $\g_p(G)\geq p+1$ and there exists an edge
$xy\in E(G^c)$ such that $\g_p(G+xy)=\g_p(G)-1$. Let $G'=G+xy$ and
$X$ be a $\g_p$-set of $G'$. Without loss of generality, assume that
$x\in X$ and $y\in \overline X$. Clearly, $y$ is a $p$-private
neighbor of $x$ with respect to $X$ in $G$ and $X\cup \{y\}$ is a
$\g_p$-set of $G$, which implies
 $$
 \mbox{$PN_p(y,X\cup \{y\},G)=\emptyset$ and $p-|N_G(y)\cap (X\cup \{y\})|=1$,}
 $$
that is, $\mu_p(y,X\cup
\{y\},G)=1$. It follows that
 $$
 r_p(G)\leq \mu_p(G)\leq \mu_p(X\cup \{y\},G)\leq \mu_p(y,X\cup
 \{y\},G)=1.$$
 Thus, $r_p(G)=\mu_p(G)=1$.
The theorem follows.\end{pf}

Note that  $|PN_p(x,X,G)|\leq deg_G(x)$ for any
$X\subseteq V(G)$ and $x\in X$. By Theorem \ref{thm4.1}, we obtain the following
corollary immediately.

\begin{cor}\label{cor4.1}
For any graph $G$ with maximum degree $\Delta(G)$ and positive
integer $p$, $r_p(G)\leq \Delta(G)+p$.
\end{cor}

\begin{cor}\label{cor4.2}
Let $p$ be a positive integer and $G$ be a graph with minimum degree
$\delta(G)$. If $\delta(G)< p$, then $r_p(G)\leq \delta(G)+p$.
\end{cor}

\begin{pf}
Let $X$ be a $\g_p$-set of $G$ and $x\in V(G)$ with degree
$\delta(G)$. Since $deg_G(x)=\delta(G)<p$, $x\in X$ by Observation~\ref{obs1.2}.
Note that $|PN_p(x,X,G)|\leq
deg_G(x)=\delta(G)$ and $p-|N_G(x)\cap X|\leq p$. By Theorem
\ref{thm4.1},
 \begin{eqnarray*}
 r_p(G)&\leq& \mu_p(G)\\
       &\leq& \mu_p(x,X,G)\\
       &=&|PN_p(x,X,G)|+\max\{0,p-|N_G(x)\cap X|\}\\
       &\leq& \delta(G)+p.
 \end{eqnarray*}
 The corollary follows.\end{pf}

Consider $p=1$. Let $X\subseteq V(G)$ and $x\in X$. If $x$ is not an isolated vertex of the induced subgraph $G[X]$,
then  $PN(x,X,G)$ defined in
(\ref{e4.1}) does not contain $x$ and $\max\{0,1-|N_G(x)\cap X|\}=0$ in (\ref{e4.4}). Otherwise, $PN(x,X,G)$ contains $x$ and
$\max\{0,1-|N_G(x)\cap X|\}=1$.
Notice that $PN_1(x,X,G)$ defined in
(\ref{e4.3}) does not contain $x$. Hence,  by (\ref{e4.5}),
 $$
 \mu_1(x,X,G)=PN_1(x,X,G)+\max\{0, 1-|N_G(x)\cap X|\}=|PN(x,X,G)|.
 $$
This fact means that
$\mu(G)$ defined in (\ref{e4.2}) is a special case of $p=1$ in
(\ref{e4.6}), that is, $\mu_1(G)=\mu(G)$. Thus, by Theorem~\ref{thm4.1}, the following
corollary holds immediately.

\begin{cor}\label{cor4.3}\textnormal{(Kok and Mynhardt~\cite{km90})}
For any graph $G$ with $\g(G)\geq 2$, $r(G)\leq\mu(G)$, with
equality if $r(G)=1$.
\end{cor}




\begin{thebibliography}{99}



\bibitem{bcf05}
M. Blidia and M. Chellali, O. Favaron, Independence and 2-domination
in trees. \emph{Austral. J. Combin.} 33 (2005) 317-327.


\bibitem{bcv06}
M. Blidia, M. Chellali and L. Volkmann, Some bounds on the
$p$-domination number in trees. \emph{Discrete Math.} 306 (2006) 2031-2037.

\bibitem{bgh08}
J.R.S. Blair, W. Goddard, S.T. Hedetniemi, S. Horton, P. Jones and
G. Kubicki, On domination and reinforcement numbers in trees.
\emph{Discrte Math.} 308 (2008) 1165-1175.



\bibitem{cfhv11}
M. Chellali, O. Favaron, A. Hansberg and L. Volkmann, $k$-domination
and $k$-independence in graphs: A survey. \emph{Graphs $\&$ Combin.}
doi 10.1007/s00373-011-1040-3.


\bibitem{cr90}
Y. Caro and Y. Roditty, A note on the $k$-domination number of a
graph, Internat. \emph{J. Math. Sci.} 13 (1990) 205-206.

\bibitem{csm03}
X. Chen, L. Sun and D. Ma, Bondage and reinforcement number of
$\g_f$ for complete multipartite graph, J. Beijin Inst. Technol. 12
(2003) 89-91.

\bibitem{dhtv98}
J. E. Dunbar, T. W. Haynes, U. Teschner and L. Volkmann, Bondage,
insensitivity, and reinforcement. Domination in Graphs: Advanced
Topics (T. W. Haynes, S. T. Hedetniemi, P. J. Slater eds.), 471-489,
Monogr. Textbooks Pure Appl. Math., 209, Marcel Dekker, New York,
(1998).

\bibitem{dl97}
G.S. Domke and R.C. Laskar, The bondage and reinforcement numbers of
$\g_f$ for some graphs. \emph{Discrete Math.} 167/168 (1997)
249-259.

\bibitem{df95}
R.G. Downey, M.R. Fellows, Fixed-parameter tractability and
completeness I: Basic results. SIAM J. Comput. 24 (1995), 873-921.

\bibitem{df97}
R.G. Downey, M.R. Fellows, Fixed-parameter tractability and
completeness II: On completeness for $W[1]$.  Theoretical Computer
Science,  54 (3) (1997), 465-474.

\bibitem{f85}
O. Favaron, On a conjecture of Fink and Jacobson concerning
$k$-domination and $k$-dependence. \emph{J. Combin. Theory Ser. B}
39 (1985) 101-102.

\bibitem{fj85}
J. F. Fink and M. S. Jacobson, $n$-domination in graphs. Graph
Theory with Applications to Algorithms and Computer Science (Y.
Alavi, A. J. Schwenk eds), 283-300, Wiley, New York, (1985).

\bibitem{gj79}
M.R. Garey and D.S. Johnson, Computers and Intractability: A Guide
to the Theory of NP-Completeness, Freeman, San Francisco, (1979).

\bibitem{hs981}
T. W. Haynes, S. T. Hedetniemi and P. J. Slater, Fundamentals of
Domination in Graphs, New York, Marcel Deliker, (1998).

\bibitem{hs982}
T. W. Haynes, S. T. Hedetniemi and P. J. Slater, Domination in
Graphs: Advanced Topics, New York, Marcel Deliker (1998).

\bibitem{hrr11}
M.A. Henning, N.J. Rad and J. Raczek, A note on total reinforcement
in graph. \emph{Discrete Appl. Math.} 159 (2011) 1443-1446.

\bibitem{hx10}
F.-T. Hu and J.-M. Xu, On the Complexity of the Bondage and
Reinforcement Problems. Journal of Complexity (2011),
doi:10.1016/j.jco.2011.11.001.


\bibitem{hwx09}
J. Huang, J.W. Wang and J.-M. Xu, Reinforcement number of digraphs.
\emph{Discrete Appl. Math.} 157 (2009) 1938-1946.


\bibitem{km90}
J. Kok and C.M. Mynhardt, Reinforcement in graphs. \emph{Congr.
Numer.} 79 (1990) 225-231.


\bibitem{ses07}
N. Sridharan, M.D. Elias and V.S.A. Subramanian, Total reinforcement
number of a graph. \emph{AKCE Int. J. Graph Comb.} 4 (2) (2007)
192-202.

\bibitem{x03}
J.-M. Xu, Theory and Application of Graphs. Kluwer Academic
Publishers, Dordrecht/Boston/London, 2003.

\bibitem{zls03}
J.H. Zhang, H.L. Liu and L. Sun, Independence bondage and
reinforcement number of some graphs. Trans. Beijin Inst. Technol.
23 (2003) 140-142.

\end{thebibliography}
\end{document}